\documentstyle[12pt]{article}
\input tcilatex

\begin{document}

\title{SEMILOGICS, QUASILOGICS AND OTHER QUANTUM STRUCTURES}
\author{V P Belavkin \\
%EndAName
The University of Nottingham, \\
Mathematics Department, NG7 2RD, UK}
\date{vpb@maths.nott.ac.uk\\
Translated from Russian by A. Pulmanova for\\
Tatra Mountains Math. Publ. {\bf 10} (1997) 199--223. Original is published
in: Mathematical Foundations of the General Systems Theory, MIEM, Moscow,
1987.}
\maketitle

\begin{abstract}
We give an axiomatic formulation of quantum structures like semilogics and
quasilogics which generalize the boolean semirings of events and fuzzy
logics. The notions of distributions, states, representations observables
and semiobservables are introduced and their Hilbert space realizations are
found. The closed and open structures in semilogics are introduced and the
regular distributions on the semilogics are studied.
\end{abstract}

\section{The logic of the systems}

\subsection{Logical relations}

In this section we summarize the logical notions, which were discussed in
connection with the foundations of quantum mechanics in [2]. The set $\ell $
of all propositions $a,b,c,\dots $ about an arbitrary system is a partially
ordered set with respect to the relation of logical inclusion: the
proposition $a\in \ell $ is included in the proposition $b\in \ell $: $a\leq
b$, if from $a$ is true always follows $b$ is true (it can be also said that
the statement $a$ is not weaker than the statement $b$ ). The least upper
bound and the greatest lower bound $\sup \{a,b\}$ and $\inf \{a,b\}$
according to the logical ordering $\leq $ are called the {\em disjunction} $%
a\vee b$ and the {\em conjunction} $a\wedge b$ of the propositions $a$ and $%
b $, respectively.

The identically false proposition, the zero 0, is the least element of $\ell$%
, i.e. $0 \leq a \quad \forall a \in \ell$, and the identically true
proposition, the unity 1, is the greatest element of $\ell$, i.e. $1 \geq a
\quad \forall a \in \ell$. The elements $0$ and $1$ are to be regarded as
nullary operations, giving the strongest and the weakest proposition in $%
\ell $.

In each logic $\ell$, the unary operation of {\em negation} $a \mapsto \bar
a $ is introduced, 
\begin{equation}
a\vee \bar a = 1 ,\qquad a\wedge \bar a = 0  \label{eq:1.1}
\end{equation}
which determines an involutive anti-automorphism of the partially ordered
set $\ell$, i.e. 
\begin{equation}
a\leq b \Rightarrow \bar b \leq \bar a , \qquad \bar{\bar a} = a
\label{eq:1.2}
\end{equation}
It is easy to prove that from $\bar{ 1} = 0$, $\bar{ 0} = 1$, the de Morgan
laws follow: 
\begin{equation}
\overline{a \vee b} = \bar a\wedge \bar b ,\qquad \overline{a \wedge b}=
\bar a \vee \bar b  \label{eq:1.3}
\end{equation}

The elements $a$, $b$, if $\bar a \geq b$ ($\bar a \leq b$), are called {\em %
disjoint} ({\em conjoint}) and the notation $a\perp b$ ($a \top b$) is used.
The disjunction of disjoint elements is denoted by a sum, and the
conjunction of conjoint elements is denoted by a product: 
\begin{eqnarray}
a\vee b = a + b & \mbox{if} & a \perp b  \label{eq:1.4} \\
a \wedge b = ab & \mbox{if} & a \top b  \nonumber
\end{eqnarray}

Note that from $a \perp b$ always follows $a\wedge b = 0$, and from $a \top
b $ follows that $a \vee b = 1$, because the inequalities $\bar a \geq b$, $%
\bar a \leq b$ give $0 = a \wedge \bar a \geq a \wedge b$, $1 = a \vee \bar
a \leq a \vee b$, correspondingly.

In classical logics, described by the boolean algebras, the distributive law 
\begin{equation}
a \wedge (b \vee c) = (a \wedge b)\vee (a\wedge c),\quad a\vee (b\wedge c)=
(a\vee b)\wedge (a\vee c),  \label{eq:5}
\end{equation}
taken as an axiom, turns every logic $\ell$, determined by the axioms~(\ref%
{eq:1.1}--\ref{eq:1.4}), in a distributive structure, where from the fact
that $a$ and $b$ are non-simultaneous, i. e. $a\wedge b = 0$, follows that $%
a $ and $b$ are mutually exclusive, i. e. $a \perp b$, and vice versa.
However, the empirical data of quantum physics show the existence of
non-simultaneous propositions $a$, $b$, that are not mutually exclusive in
the sense that $\bar a \geq b$. For example, the statements ``the impulse of
a particle has the value $p \in \Delta p$'' and ``the value of the
coordinate is $q \in \Delta q$'' are, according to the uncertainty law,
non-simultaneous, though they don't exclude each other: if the result of
some experiment would be that one of the statements is true, it doesn't mean
that the other is false at all. Instead of the completely non-obvious
condition (5), in the logics of physical systems a much weaker and
intuitively clear condition holds: 
\begin{equation}
(a\vee b)\wedge c = a\vee (b\wedge c), \qquad \forall a\leq \bar b \leq c
\label{eq:1.5}
\end{equation}
Its meaning will be revealed in the following theorem.

\begin{theorem}
Each segment $[a,c]$ in a logic $\ell $, satisfying the condition~(\ref
{eq:1.5}), is a logic, where the zero equals $a$, the unity equals $c$ and
the operation of complementation is given by 
\[
b\mapsto a+(\bar{b}c)=(a+\bar{b})c 
\]
\label{th:2.1}
\end{theorem}

\noindent{\sc Proof} First, with respect to~(\ref{eq:1.5}), the
complementation is well defined and the equalities $b \vee (a\vee \bar b
\wedge c) = c$ and $b \wedge (a\vee \bar b \wedge c)=a$ follow straight from
the definitions of supremum and infimum $\vee$ and $\wedge$. As $a \vee \bar
b \wedge c$ is a monotone function of $b$, the inequality $a\leq b_1\leq b_2
\leq c$ leads to the inequality $a\vee \bar b_1 \wedge c\geq a\vee \bar b_2
\wedge c$.

Now it stays to prove the involutivity 
\[
a \vee \overline{(a\vee \bar b\wedge c)}\wedge c = a\vee \bar a \wedge b
\vee \bar c \wedge c = b 
\]
The proof is complete.

Thus the system of logical axioms~(\ref{eq:1.1}--\ref{eq:1.5}) is natural,
but not independent, as follows from the de Morgan laws~(\ref{eq:1.3}).

\subsection{Semilogics and quasilogics}

A {\em semilogic} is a mathematical model ${\cal M}$ with the zero $0\in 
{\cal M}$, a partial operation of the logical multiplication $(a,b)\mapsto
ab\in {\cal M}$ and an operation of the logical addition $\{a_i\}\mapsto
\sum_ia_i\in {\cal M}$, defined for some families $\{a_i\}$ of pairwise
orthogonal elements $a_i\in {\cal M}$, $a_{i_1}a_{i_2}=0$, $i_1\neq i_2$;
described by the following axioms: 
\[
0a=0,\ aa=a,\ ab=ba,\ (ab)c=(bc)a,\ \sum_iaa_i=a\sum_ia_i 
\]
The binary relation $a\bar{\downarrow}b\Leftrightarrow ab\in {\cal M}$,
which can be seen as the domain of the commutative multiplication and is
called the relation of logical commutation, has the following properties: 
\begin{eqnarray*}
&\{a\in {\cal M}:0\bar{\downarrow}a\}={\cal M},& \\
&a\bar{\downarrow}a,\quad \forall a\in {\cal M}\qquad &\mbox{(reflexivity)}
\\
&a\bar{\downarrow}b=b\bar{\downarrow}a,\qquad &\mbox{(symmetry)} \\
&a\bar{\downarrow}b,\ c\bar{\downarrow}a,\ b\bar{\downarrow}c\Rightarrow ab%
\bar{\downarrow}c\qquad &\mbox{(associativity)}
\end{eqnarray*}
The summable families $\{a_i\}$, for which the symmetric and associative
sums $\sum_ia_i$ are defined, consist of either mutually different, or zero
elements $a_i\in {\cal M}$ and can be identified with the subsets $%
A\subseteq {\cal M}$, $0\notin A$, where $\sum_ia_i=\sum_{a\in A}a$. The
family ${\cal A}$ of such subsets contains the empty set $\emptyset \subset 
{\cal M}$, here, according to the definition, $\sum_{a\in \emptyset }a=0$;
all the nonzero one-point subsets $\{a\}\subset {\cal M}$, $a\neq 0$, here $%
\sum_{a\in \{a\}}a=a$; while the union of arbitrary summable subsets $A_i\in 
{\cal A}$, determining summable elements $a_i=\sum_{a\in A_i}a$, $%
\{a_i\}=A\in {\cal A}$, is also summable, i. e. $\cup A_i\in {\cal A}$.
Besides, it is assumed that if $a\bar{\downarrow}b$, then decompositions $%
a=\sum_{a^{\prime }\in A}a^{\prime }$ and $b=\sum_{b^{\prime }\in
B}b^{\prime }$, included in some decomposition $c=\sum_{c^{\prime }\in
C}c^{\prime }\in {\cal M}$, $A\cup B\subseteq C$ exist, the intersection of
which gives the decomposition $ab=\sum_{d\in A\cap B}d$. In particular, this
means that the logical multiplication and addition define one and the same
relation of logical ordering: $ab=a$ iff for any decomposition $%
a=\sum_{a^{\prime }\in A}a^{\prime }$, $A\in {\cal A}$, there is a
decomposition $b=\sum_{b^{\prime }\in B}b^{\prime }$, $B\in {\cal A}$, such
that $A\subseteq B$, denoted by $a\leq b$. By the axioms, the binary
relation $\leq $ is reflexive, transitive and antisymmetric, $0$ is the
least element and the product $ab$ is equal to the greatest lower bound $%
a\wedge b=\inf \{a,b\}$ and the sum $\sum_{a\in A}a$ is equal to the least
upper bound $\bigvee_{a\in A}a=\sup_{a\in A}a$ with respect to this ordering.

It is usual to assume that $a\wedge b\in {\cal M}$ for any $a,b\in{\cal M}$,
and only the finite sums $\sum_i a_i$, determined by finite subsets $A
\subset {\cal M}$ are studied, but it is possible to consider infinite sums
too, defining those as least upper bounds $\sum_i a_i = \bigvee_i a_i$.

An element $a\in {\cal M}$ is called an {\em atom}, if from $b<a$ follows $%
b=0$, i. e. if this element cannot be decomposed into a sum of two nonzero
elements. If for any $b\neq 0$ there is an atom $a\leq b$, the semilogic is 
{\em atomic}. Each semilogic of finite rank is atomic, as any element $a\in 
{\cal M}$ can be represented in a form of a sum $\sum_ia_i=a$ of atoms $%
a_i\in {\cal M}$. A semilogic, such that for any family $\{a_i\}_{i\in I}$
of cardinality $|I|\leq \tau $ the greatest lower bound $\bigwedge_{i\in
I}a_i\in {\cal M}$ is defined, is called $\tau $-{\em complete}, or a $\tau $%
-{\em semilogic} if $\tau $ is not less than the maximum cardinality $|A|$
of summable families $A\in {\cal M}$.

A {\em quasilogic}\footnote{%
The quasilogics are also called D-posets, see for example [3].} is a
partially ordered set ${\cal N}$ with a partial operation of subtraction: 
\begin{equation}
a\leq b\Rightarrow \exists b-a\leq b,\quad b-\left( b-a\right) =a
\label{ax:0}
\end{equation}
having the property 
\begin{eqnarray}
b\leq c &\Rightarrow &b-a\leq c-a,(c-a)-(b-a)=c-b,\quad \forall a\leq b
\label{ax:1} \\
a\leq b &\Rightarrow &c-b\leq c-a,(c-a)-(c-b)=b-a,\quad \forall c\geq b
\label{ax:2}
\end{eqnarray}
If the greatest element $1\geq a$, $\forall a\in {\cal N}$ exists in ${\cal N%
}$, then we can define the complementation $\bar{a}=1-a$, which obviously
has the property~(\ref{eq:1.2}); however, this is not necessary in general
and we will suppose just that ${\cal N}$ is upwards directed in the sense of
existence of a majorating element $c\in {\cal N}$ : $c\geq a,b$ for each
pair $a,b\in {\cal N}$. This yields the existence and uniqueness of the
least element $0=a-a$ and an involutive anti-automorphism $a\mapsto c-a$ of
any segment $[0,c]\subset {\cal N}$ onto itself, so that if the existence of
supremum and infimum $a\vee b$ and $a\wedge b$ is assumed, the de Morgan
formulas hold: 
\begin{equation}
c-a\vee b=(c-a)\wedge (c-b),\quad c-a\wedge b=(c-a)\vee (c-b),\quad \forall
c\geq a,b  \label{eq:1.6}
\end{equation}
For each pair $a,b\in {\cal N}$, such that there is an element $c\geq a,b$ : 
$c-a\geq b$ (or, equivalently, $c-b\geq a$, denoted by $a\underline{\uparrow 
}b$) this supposition allows us to define a partial operation of addition $%
a+b$ by putting 
\begin{equation}
c-((c-a)-b)=a+b=c-((c-b)-a),\quad \forall a\underline{\uparrow }b
\label{eq:1.7}
\end{equation}
First, $(a+b)-b=a$ and $(a+b)-a=b$ hold by the axioms~(\ref{ax:1}) and~(\ref%
{ax:2}) and the equality $(c-a)-b=(c-b)-a$ follows. Moreover, the definition
of the sum $a+b$ doesn't depend on $c$, because for any $c_1$ and $c_2$
there is an upper bound $c\geq c_1,c_2$, such that by the axioms~(\ref{ax:1}%
) and~(\ref{ax:2}) 
\[
c_1-((c_1-a)-b)=c-((c-a)-b)=c_2-((c_2-a)-b) 
\]

It is easy to prove that such binary partial operation of addition is
commutative: if $a\underline{\uparrow }b$ then $b\underline{\uparrow }a$ and 
$a+b=b+a$, associative: if $a\underline{\uparrow }(b+c)$ then $(a+b)%
\underline{\uparrow }c$ and 
\[
a+(b+c)=(a+b)+c, 
\]
zero $0$ is the neutral element $0+a=a+0=a$ and 
\begin{equation}
a+b=a\vee b+a\wedge b\qquad \mbox{if}\ a\vee b,a\wedge b\in {\cal N}
\label{eq:1.8}
\end{equation}

We will say that $a,b\in {\cal N}$ {\em quasicommute}, $a\downarrow b$, if
there is a $c\geq a,b$ such that $c-a\leq b$ (or $c-b \leq a$). The relation 
$a\downarrow b$ means the summability $p\underline{\uparrow}q$ of the
elements $p=c-a$ and $q=c-b$, determining the decompositions $a=(ab)_c+q$
and $b=(ab)_c+p$ where $(ab)_c$ denotes the {\em quasiproduct of $a$ and $b$
with respect to $c$}, defined as 
\[
a-(c-b)=(ab)_c=b-(c-a) \qquad \forall a\downarrow b 
\]

The elements $a,b\in {\cal N}$ are {\em disjoint} $a\perp b$, if they are
summable: $a\underline{\uparrow}b$ and non-simultaneous: $a\wedge b=0$; and 
{\em (logically) commuting} $a\bar{\downarrow}b$, if for some $c\geq a,b$
the elements $p=c-a$ and $q=c-b$ are disjoint; in agreement with~(\ref%
{eq:1.8}) we have $a+b=a\vee b$ if $a\perp b$ and $ab=a\wedge b$ if $a\bar{%
\downarrow}b$. A quasilogic, where only the disjoint elements are summable
is a semilogic with respect to the defined partial operations of logical sum 
$a+b=a\vee b$ and product $ab=a\wedge b$ and is called an {\em ortho-logic}
or simply a {\em logic}.

Hence a logic $\ell$ is a upwards directed semilogic ${\cal M}$ with an
operation of subtraction $b-a\in \ell$ if $ab=a$, behaving by the axioms~(%
\ref{ax:1}) and (\ref{ax:2}), where the summable families are at most all
finite families $\{ a_i\}$ of mutually disjoint elements $a_{i_1}\perp
a_{i_2}$, $i_1\neq i_2$: 
\[
\sum_{i=1}^na_i = \sum_{i=1}^{n-1}a_i +a_n\qquad \forall n<\infty 
\]

A quasilogic (logic) ${\cal N}$ is $\tau$-{\em complete}, if each downwards
directed family $\{a_i\}_{i\in I}$ of cardinality $\tau$ has the greatest
lower bound $\bigwedge a_i$, and a $\tau$-{\em quasilogic} if for each
family $\{ a_i\}_{i\in I}$, $|I|\leq\tau$ of mutually summable elements $a_i
\in {\cal N}$ there is an element $a\in {\cal N}$: $a_i\leq a$, $\forall
i\in I$.

\subsection{Homomorphisms and closures}

Functional relations between different systems or their parts are described
by {\em logical homomorphisms} of the corresponding logics, defined as
additive mappings $h:\ {\cal M}\rightarrow {\cal N}$, 
\begin{equation}
h(\sum a_i)=\sum h(a_i)  \label{eq:1.9}
\end{equation}
giving a representation of the logic (semilogic) ${\cal M}$ in a logic
(quasilogic) ${\cal N}$. Each such mapping is monotone: 
\[
a\leq b\Rightarrow h(a)\leq h(b), 
\]
preserves zero: $h(0)=0$ and the operation of subtraction: 
\[
h(b-a)=h(b)-h(a) 
\]
if defined in ${\cal M}$, and also the relation of commutation 
\[
a\downarrow b\Rightarrow h(a)\downarrow h(b) 
\]
if ${\cal N}$ is a logic; in this case $h(ab)=h(a)h(b)$.

A homomorphism $h$ is $\tau $-{\em additive (completely additive)}, if (13)
is fulfilled not only for all finite summable families, but also for all
infinite summable families $\{a_i\}\subset {\cal M}$ of cardinality $|I|\leq
\tau $ (of arbitrary cardinality); and a $\tau $-{\em homomorphism (normal
homomorphism)}, if $\bigwedge_{i\in I}h(a_i)=0$ for $\bigwedge_{i\in I}a_i=0$
for each zero directed family $\{a_i\}$ of cardinality $|I|\leq \tau $ (any
cardinality). Each $\tau $-homomorphism of a $\tau $-logic ${\cal M}$ into a 
$\tau $-logic ${\cal N}$ is $\tau $-additive.

Note that in the definition of the homomorphisms we, in general, do not need
preservation of unity, if it occurs in ${\cal M}$ and ${\cal N}$; however,
if $h$ is an epimorphism, i. e. if $h({\cal M})={\cal N}$, then from $1\in 
{\cal M}$ follows that $h(1)=1\in {\cal N}$. An {\em approximate unity} in a
semilogic (quasilogic) ${\cal M}$ is a family ${\cal I}\subseteq {\cal M}$
of elements $i\in {\cal I}$, satisfying the following conditions:

\begin{enumerate}
\item for any $a\in {\cal M}$ there is an element $i\in {\cal I}$, $i\geq a$,

\item for each pair $i_1$, $i_2$ there is an $i\geq i_1,i_2$.
\end{enumerate}

A homomorphism $h:{\cal M}\rightarrow {\cal N}$ is called {\em approximate},
if the least logic (quasilogic) $\ell \subseteq {\cal N}$, containing the
image $h({\cal M})$, is an approximate unity in ${\cal N}$; and $\tau $-{\em %
approximate} if the least $\tau $-logic ($\tau $-quasilogic) $\ell \subseteq 
{\cal N}$, $h({\cal M})\subseteq \ell $ is an approximate unity. If ${\cal M}
$ is a $\tau $-logic and $h$ is a $\tau $-homomorphism, then $h({\cal M})$
is a $\tau $-quasilogic and the notions of approximity and $\tau $%
-approximity are identical; here $h(1)=1\in {\cal N}$, if $1\in {\cal M}$.
Any homomorphism $h:{\cal M}\rightarrow {\cal N}$, preserving the
approximate unity ${\cal I}\subseteq {\cal M}$: $h({\cal I})={\cal J}$,
where ${\cal J}\subset {\cal N}$ is a given approximate unity, is
approximate.

A subset ${\cal I}\subseteq {\cal M}$ is called an {\em upper family} of a
semilogic (quasilogic) ${\cal M}$, if an element $i\geq a$ exists in ${\cal I%
}$ for each element $a\in {\cal M}$ and for any $i_1$, $i_2\in {\cal I}$, $%
i_1,i_2\geq a$, there is $i\geq a$ in ${\cal I}$, such that $i_1,i_2\geq i$.
The {\em lower family} ${\cal K}\subseteq {\cal M}$ is defined analogically,
such family is called {\em opposite} to ${\cal I}$, if $i-k\in {\cal I}$
whenever $i\in {\cal I}$, $k\in {\cal K}$, $i\geq k$. These families ${\cal I%
}$ and ${\cal K}$ are called {\em approximate} from above and from below,
respectively, if 
\[
\bigwedge \{i\in {\cal I}\ :\ i\geq a\}=a,\bigvee \{k\in {\cal K}\ :\ k\leq
a\}=a 
\]
for $a\in {\cal M}$; and $\tau $-{\em approximate}, if the above conditions
are fulfilled for some subsets $\{i\}\geq a$, $\{k\}\leq a$ of ${\cal I}$
and ${\cal K}$ with cardinality $\leq \tau $.

A {\em closure} in a semilogic (logic) is a projection $k:\ {\cal M}%
\rightarrow {\cal M}$, $k\circ k=k$, satisfying the conditions 
\[
k(0)=0,\quad k(a)\geq a,\quad k(a\vee b)=k(a)\vee k(b). 
\]
An element $a\in {\cal M}$ is {\em closed}, if it equals its closure: $%
a=k(a) $, and {\em open}, if $k-a$ is closed whenever $k\geq a$ is closed.
The set ${\cal K}=k({\cal M})$ of all closed elements is a lower family,
containing the supremum $k_1\vee k_2$ of any pair $k_1,k_2\in {\cal K}$ and
the infimum $\bigwedge k_i$ of any family $\{k_i\}\subset {\cal K}$ if it is
defined in ${\cal M}$; and is $\tau $-approximate from below, if some upper $%
\tau $-approximate family of open elements $i\in {\cal M}$ is contained in $%
{\cal M}$. The set ${\cal I}$ of all open elements $i\in {\cal M}$, defining
a projection 
\[
i(a)=\bigvee \{i\in {\cal I}\ :\ i\leq a\}, 
\]
fulfils the condition $i_1\wedge i_2\in {\cal I}$ for each pair $i_1,i_2\in 
{\cal I}$ and is an upper family; if it forms an approximate unity in ${\cal %
M}$, then it is called {\em dual} to ${\cal K}$. If ${\cal M}$ contains an
unity $1\in {\cal M}$, then $1\in {\cal I}$, because the family dual to $%
{\cal K}$ is a set of the complements $i=1-k$, $k\in {\cal K}$.

\subsection{Logical states}

A real positive function $m$ on a semilogic (quasilogic) ${\cal M}$ is
called a {\em distribution}, if it satisfies the additivity condition 
\begin{equation}
m(\sum a_i)=\sum m(a_i)  \label{eq:1.II}
\end{equation}
Any distribution is monotone: $a\leq b\Rightarrow m(a)\leq m(b)$, preserves
zero: $m(0)=0$ and can be considered as a homomorphism into the quasilogic $%
R_{+}$ of positive numbers with the natural ordering and with the operations
of subtraction and addition. The distributions defined on a $\tau $%
-semilogic ($\tau $-quasilogic) are called $\tau $-{\em additive} ($\tau $-%
{\em distributions}), if they are $\tau $-additive homomorphisms ($\tau $%
-homomorphisms); and {\em measures} ($\tau $-{\em additive measures, $\tau $%
-measures}), if ${\cal M}$ is a logic ($\tau $-logic) $\ell $. A
distribution (measure) $m:{\cal M}\rightarrow R_{+}$ is {\em bounded}, if it
has a finite mass $\Vert m\Vert =\sup \sum m(a_i)$ (the supremum is taken
over all possible disjoint families $\{a_i\}\subset {\cal M}$); and {\em %
completely additive} ({\em normal}) if it is $\tau $-additive (a $\tau $%
-distribution) for any $\tau $. Each distribution defined on a logic
(quasilogic) with unity $1\in {\cal M}$ is finite: $\Vert m\Vert
=m(1)<\infty $, and completely additive, if it is $\tau $-additive for some $%
\tau =\dim ({\cal M})$, where $\dim ({\cal M})$ is the maximal number of
pairwise summable elements in ${\cal M}$.

The set of all elements $a\in {\cal M}$, such that $m(a)=0$, forms an ideal
in ${\cal M}$, called the {\em zero ideal} of the distribution (measure) $m$%
. In general, an {\em ideal} of a semilogic ${\cal M}$ is a nonempty proper
subset $\Delta \subset {\cal M}$, having the properties:

\begin{description}
\item[a)] $a\wedge b\in \Delta $ for each $a\in {\cal M}$ and $b\in \Delta $,

\item[b)] $\sum a_i\in \Delta $ for all summable $a_i\in \Delta $.
\end{description}

The ideal $\Delta $ is called a $\tau $-{\em ideal} (a {\em complete ideal}%
), if $\sum a_i\in \Delta $ for all summable families $\{a_i\}_{i\in I}$ of
cardinality $|I|\leq \tau $ (any cardinality). The zero ideal of a $\tau $%
-additive (completely additive) distribution is a $\tau $-ideal (a complete
ideal). Other examples of $\tau $-ideals in semilogics of range $\tau $ ($%
\tau $-complete logics) are the {\em regular ideals}, determined by some
(not maximal) element $b\in {\cal M}$ as the subset $\Delta =\{a\in {\cal M}%
\ :\ a\perp b\}$. The ideal $\Delta $ is called a {\em maximal ideal}, if it
is not a proper subset of any other ideal. Each maximal regular ideal is
determined by some atom $q\in {\cal M}$.

A distribution $p:\ {\cal M}\rightarrow [0,1]$ of a mass $\Vert p\Vert =1$
is a {\em probability distribution (measure)}. A probability distribution
(measure) on a logic (quasilogic) is a {\em state} if $p(c)=1$ for some $%
c\in {\cal M}$. The set $\nabla $ of all such $c$, the {\em support} of $p$,
doesn't contain zero $0\notin \nabla $ and has the following properties:

\begin{description}
\item[a)] from $b\in \nabla $, $a\geq b$ follows that $a\in \nabla $,

\item[b)] from $a,b\in \nabla $ and $c-a\underline{\uparrow }c-b$ for some $%
c\geq a,b\in \nabla $ follows that $c-\left[ (c-a)+(c-b)\right] \in \nabla $.
\end{description}

In the case when ${\cal M}$ is a semilogic, the last condition has the form 
\[
a,b\in \nabla ,a\bar{\downarrow}b\Rightarrow ab\in \nabla , 
\]
where $\bar{\downarrow}$ is the relation of logical commutation and $%
ab=a\wedge b$ is the logical multiplication with the domain $\bar{\downarrow}
$. Such set $\nabla $ in a logic $\ell $ is called a {\em filter}, and a $%
\tau $-{\em filter} ({\em complete filter} if, together with the condition
b), the condition $\bigwedge_{i\in I}a_i\in \nabla $ is satisfied for any
pairwise commuting family $\{a_i\}_{i\in I}$ of cardinality $|I|\leq \tau $
(any cardinality). The unit filters - the supports of a $\tau $-state $p$,
as well as the regular filters in a $\tau $-logic $\ell $, defined as
subsets $\nabla =\{a\in \ell \ :\ a\geq b\}$ for some $b\in \ell $, are $%
\tau $-filters.

Let ${\cal M}$ be a semilogic (logic) with an approximate upper family $%
{\cal I}$ and an approximate lower family ${\cal K}$. A distribution
(measure, state) is {\em regular} if it is continuous with respect to the
topology on ${\cal M}$, generated by the intervals $[k,i]=\{a\ :\ k\leq
a\leq i\}$ for all $i\in {\cal I}$, $k\in {\cal K}$. A necessary and
sufficient condition for regularity of a distribution is given by the
following approximate property: 
\[
\sup \{m(k)\ :\ k\leq a,\ k\in {\cal K}\}=m(a)=\inf \{m(i)\ :\ i\geq a,\
i\in {\cal I}\}, 
\]
in the case that ${\cal K}$ is opposite to ${\cal I}$, the second part of
the above is a consequence of the first part. A distribution $m$ is called $%
\tau $-{\em regular}, if $m(a_\lambda )\rightarrow m(a)$ for any $\tau $-net 
$\{a_\lambda \}$ converging to $a\in {\cal M}$, in particular, a
distribution is $\delta $-regular if $m(a_n)\rightarrow m(a)$ for any $%
a_n\rightarrow a$. Each $\tau $-distribution on a $\tau $-complete semilogic 
${\cal M}$ is $\tau $-regular, and completely regular, if ${\cal I}$ and $%
{\cal K}$ are $\tau $-approximate families in ${\cal M}$. In particular, any
normal distribution is regular, the converse holds for ${\cal I}={\cal M}=%
{\cal K}$.

\section{Classical systems}

\subsection{Boolean rings}

A quasilogic ${\cal N}$ with a trivial relation of quasicommutation $%
a\downarrow b$, $\forall a,b\in {\cal N}$, is called a {\em quasiring}, and
a {\em (boolean) ring}, if each pair $a,b\in {\cal N}$ logically commute, i.
e. $a\bar{\downarrow}b$. In other words, a quasiring (ring) is a partially
ordered set ${\cal N}$ with subtraction: $a\leq b\Rightarrow b-a\in {\cal N}$%
, such that for any pair $a,b$ an element $c\geq a,b$ exists, such that the
elements $c-a$ and $c-b$ are summable (disjoint): $c-a\geq b$ and $c-b\geq a$
($c-a\perp c-b$). Classical systems, ruled by the boolean logic, are usually
described by distributive logics, i.e. boolean rings ${\cal A}$ (algebras)
of propositions, taking the conjunction $a\wedge b=ab$, the difference $%
a\setminus b=a-ab$ and the disjunction $a\vee b=ab+a\setminus b+b\setminus a$
as the principal operations, defined everywhere by triviality of the
relation $\bar{\downarrow}$. A boolean ring (a {\em boolean algebra} is a
ring with unity $1\in {\cal A}$, $1\wedge a=a$, $\forall a\in {\cal A}$) is
a $\tau $-{\em complete ring} (a {\em $\tau $-ring, $\tau $-algebra}) if it
is a $\tau $-complete logic (a $\tau $-logic) with respect to the
multiplication $ab=a\wedge b$, addition $a+b=a\vee b$ if $a\wedge b=0$ and
subtraction $a-b=a\setminus b$ for $a\wedge b=b$. In a $\tau $-complete
ring, the $\tau $-infinite conjunction $\bigwedge_{i\in I}a_i\in {\cal A}$,
for arbitrary $a_i\in {\cal A}$, if $|I|\leq \tau $, and, consequently, the $%
\tau $-infinite disjunction $\bigvee_{i\in I}a_i\in {\cal A}$, if $a_i\leq
a\in {\cal A}$, $\forall i\in I$, are defined ($\bigvee_{i\in I}a_i\in {\cal %
A}$ for arbitrary $a_i\in {\cal A}$ if $|I|\leq \tau $ for the $\tau $%
-rings). By method of transfinite induction with respect to the cardinal
number $\tau $, it can be proved that the existence of the least upper bound 
$\bigvee_{i\in I}a_i=\sum_{i\in I}a_i$ for an arbitrary majorated $\tau $%
-family $\{a_i\}_{i\in I}\subset {\cal A}$, $|I|=\tau $, of mutually
disjoint elements $a_i\leq a$, $\forall i\in I$, is a sufficient condition
for $\tau $-completeness of the ring ${\cal A}$; moreover, ${\cal A}$ is a $%
\tau $-ring if the majoration condition is not substantial.

The {\em subalgebras} and {\em factor-algebras} of $\tau$-rings are the $%
\tau $-closed and the $\tau$-con\-gru\-ent subsets $({\cal A})\subseteq 
{\cal A}$ and $\{ {\cal A}\} \preceq {\cal A}$ with the natural boolean
operations. A $\tau$-{\em closed subset} of ${\cal A}$ is a subset $({\cal A}%
)$ containing the zero of the ring ${\cal A}$, $0\in {\cal A}$, which is
closed under the difference $a,b\in ({\cal A}) \Rightarrow a\setminus b\in (%
{\cal A})$ and the $\tau$-infinite conjunction $\bigwedge_{t\in T}a_t\in (%
{\cal A})$, if $a_t\in {\cal A}$, $\forall t\in T$, $|T|\leq \tau$. A $\tau$-%
{\em congruent factor-set} in ${\cal A}$ is a set of equivalence classes $\{ 
{\cal A}\}\subset {\cal A}$ of an equivalence given by some ideal $%
\Delta\subset {\cal A}$: $a\sim b$ iff $a\setminus b\in \Delta$ and $%
b\setminus a\in \Delta $.

Usually, as $\tau $-complete rings, the rings ${\cal A}$ of subsets $%
A\subseteq X$ of some set $X$ with natural boolean operations, i.e. some $%
\tau $-closed subalgebras ${\cal A}\subseteq {\cal A}(X)$ of the field $%
{\cal A}(X)$ of all subsets of $X$, not necessarily containing the unity $%
X\in {\cal A}(X)$, are considered. The points $x\in X$ of this set,
contained at least in one $A\in {\cal A}$, are interpreted as the elementary
propositions, defining complete descriptions of the states of some classical
system, resulting from attaching the true value to some of these points $x$,
so that $A\in {\cal A}$ is true if $x\in {\cal A}$, and false if $x\notin 
{\cal A}$.

A ring ${\cal A}$ is {\em separating}, if for each $x_1\neq x_2$, an $A\in 
{\cal A}$ exists, such that $x_1\in A$ and $x_2\notin A$, and {\em perfect},
if each non-trivial maximal ideal $\Delta\subset {\cal A}$ is generated by
some point $x\in X$, i.e. $\Delta = \{ A\in {\cal A}\ :\ x\notin A\}$. A $%
\tau$-complete ring is $\tau$-{\em perfect}, if the last condition is
fulfilled at least for all maximal $\tau$-ideals in ${\cal A}$. Each $\tau$%
-complete ring ${\cal A} \subseteq {\cal A}(X)$ is isomorphic with some
separating $\tau$-perfect ring ${\cal B}\subseteq{\cal A}(X)$, which is a
boolean $\tau$-algebra with unity $E\in {\cal B}$, if there is a unity in $%
{\cal A}$ (unity means a maximal subset, not necessarily identical with $X$).

The homomorphisms of boolean rings ${\cal A}\subseteq {\cal A}(X)$, ${\cal B}%
\subseteq {\cal A}(E)$, are the logical (algebraical) homomorphisms $h:\ 
{\cal B}\rightarrow {\cal A}$ of them as logics (algebras), induced by the
partial mappings $f:\ (X)\rightarrow E$, $h=f^{*}$, where $f^{*}(B)=\{x\ :\
f(x)\in B\}$ is the preimage of the subset $B\subseteq E$. Such mapping $f$
satisfy the condition $f^{*}(B)\in {\cal A}$ for each $B\in {\cal B}$ and
are called $({\cal A},{\cal B})$-measurable. A sufficient condition for an
algebraical $\tau $-homomorphism $h:\ {\cal B}\rightarrow {\cal A}$ to be
induced by some measurable mapping $f:\ (X)\rightarrow E$ is the $\tau $%
-perfectness of the ring ${\cal B}$. If also ${\cal A}$ is a $\tau $-perfect
ring, then this condition is necessary, too. Each $\tau $-subalgebra of a $%
\tau $-ring ${\cal A}\subseteq {\cal A}(X)$ is obviously a ring of subsets
of $X$, but not every factor-algebra ${\cal A}$ can be represented as a
subalgebra ${\cal A}_1\subseteq {\cal A}(X_1)$. For a factor-algebra of a $%
\tau $-ring ${\cal A}\subseteq {\cal A}(X)$ by the $\tau $-ideal $\Delta $,
a sufficient, and, in the case that this ring is $\tau $-perfect, also a
necessary condition to be a $\tau $-field is that the ideal is generated by
some subset $X_1\subset X$, i.e. $\Delta =\{A\in {\cal A}\ :\ A\cap X_1=0\}$%
. As for the $\sigma $-rings, the following theorem shows to be useful (see
[4]):

\begin{theorem}
Let ${\cal B}$ be a standard field of borel subsets $B\in {\cal B}$ of a
polish space $E$. Then each $\sigma $-homomorphism $h:\ {\cal B}\rightarrow 
{\cal A}/\Delta $ into a factor algebra of some $\sigma $-ring ${\cal A}%
\subseteq {\cal A}(X)$ by a $\sigma $-ideal $\Delta \subset {\cal A}$ is
induced by a point mapping $f:\ X\rightarrow E$ in the sense that $h=jf^{*}$%
, where $j:\ A\rightarrow \langle A\rangle $ is the factorisation
homomorphism ${\cal A}\rightarrow {\cal A}/\Delta $. \label{th:2.2}
\end{theorem}

\subsection{Semirings and topologies}

A {\em boolean semiring} ${\cal B}$ is a distributive semilogic 
\[
a\wedge (\bigvee_{i\in I}b_i)=\bigvee_{i\in I}(a\wedge b_i),\qquad a\in 
{\cal B},\quad b=\bigvee_{i\in I}b_i\in {\cal B} 
\]
with respect to the everywhere defined operation of conjunction $a\wedge
b=ab $ and the operation of disjunction $\bigvee_{i\in I}b_i$ defined for
summable families $\{b_i\}$ of pairwise non-simultaneous elements $%
b_{i_1}\wedge b_{i_2}=0$, $i_1\neq i_2$ as a logical sum: $\bigvee b_i=\sum
b_i$. Moreover, a {\em semiring of rank $\tau $} is a boolean semiring,
which is a semilogic of rank $\tau $, and a $\tau $-{\em semiring} is a
semiring of rank $\tau $, such that the greatest lower bounds $%
\bigwedge_{i\in I}b_i$ are defined for families of cardinality $|I|\leq \tau 
$.

Each semiring ($\tau$-semiring), for which the disjunction is defined for
each disjoint pair $a,b\in {\cal B}$ (for each disjoint family $\{ b_i\}$ of
cardinality $|I|\leq\tau$), is a ring (a $\tau$-ring), because the relation $%
a\leq c$ in this case means the existence of $b=c-a$, which is, by
distributivity of ${\cal B}$, uniquely defined by the conditions $a\wedge
b=0 $, $a\vee b=c$.

A semiring of subsets is any boolean semiring ${\cal B}\subseteq{\cal A}(X)$
of some set $X$ with the empty set $\emptyset\in {\cal B}$ as zero, the
intersection $A\cap B=AB\in {\cal B}$ for each $A,B\in {\cal B}$ and union $%
\bigcup B_i=\sum B_i\in {\cal B}$ for summable families $\{ B_i\}$ of
pairwise disjoint subsets $B_{i_1}\cap B_{i_2}=\emptyset$, $i_1\neq i_2$.
Each boolean semiring of subsets ${\cal B}\subseteq {\cal A}(X)$ generates a
boolean ring (a $\tau$-ring) of subsets ${\cal A}\subseteq {\cal A}(X)$,
consisting of all finite ($\tau$-infinite) sums $A=\sum B_i$, $B_i\in {\cal B%
}$, and $\tau$-complete if ${\cal B}$ is a $\tau$-complete semiring.

A {\em system of open subsets} in a semiring ${\cal B}\subseteq {\cal A}(X)$
is a family ${\cal I}\subseteq {\cal B}$ of subsets $I\subseteq X$,
satisfying the conditions:

\begin{description}
\item[1)] for each $B\in {\cal B}$ there is a majorating subset $I\in {\cal I%
}$, $B\subseteq I$,

\item[2)] if $I_1,I_2\in {\cal I}$ then $I_1\cap I_2\in {\cal I}$,

\item[3)] $I(B)=\bigcup_{I\subseteq B}I\in {\cal I}$ for each $B\in {\cal B} 
$.
\end{description}

The mapping $B\mapsto I(B)$ maps each subset $B\in {\cal B}$ to the greatest
open subset $I(B)\in {\cal I}$ contained in it, the set $I(B)$ is called the 
{\em interior} of $B$. The mapping has the following properties: $I(B_1\cap
B_2)=I(B_1)\cap I(B_2)$, $I(\cup B_i)\supseteq \cup I(B_i)$, while $I(B)=B$
iff $B\in {\cal I}$ ($B$ is open). The family ${\cal I}$ contains all unions 
$\cup I_i$ of open subsets $I_i\in {\cal I}$ defined in ${\cal B}$,
including the empty set $\emptyset \in {\cal I}$, it is an approximate unity
if ${\cal B}$ is a ring and it is called {\em approximate from above}, if,
instead of 1), a stronger condition is satisfied: $\bigcap_{I\supseteq B}I=B$
for every $B\in {\cal B}$.

A {\em system of closed subsets} in a semiring ${\cal B}\subseteq {\cal A}%
(X) $ is a family ${\cal K}\subseteq {\cal B}$ of subsets $K\subseteq X$,
fulfilling the conditions:

\begin{description}
\item[1)] $\emptyset \in {\cal K}$,

\item[2)] if $K_1,K_2\in {\cal K}$, then $K_1\cap K_2\in {\cal K}$,

\item[3)] $K(B)=\bigcap_{K\supseteq B}K\in {\cal K}$ for each $B\in {\cal B} 
$.
\end{description}

The mapping $k:\ B\mapsto k(B)$ is the closure mapping, assigning to every $%
B\in {\cal B}$ the least closed subset $K\in {\cal K}$, containing $B$. The
family ${\cal K}$ contains all intersections $\cap K_i$ of closed subsets $%
K_i\in {\cal K}$ defined in ${\cal B}$, it is an approximate unity if ${\cal %
B}$ is a ring and it is called {\em approximate from below} if in place of
1) we have $\bigcup_{K\subseteq B}K=B$ for each $B\in {\cal B}$.

A {\em topology} in a boolean ring ${\cal B}$ is a pair ${\cal T}=({\cal I},%
{\cal K})$ of families ${\cal I},{\cal K}\subseteq {\cal B}$ of open and
closed subsets $I,K\subseteq X$, which is in duality: $I\cap \bar{K}\in 
{\cal I}$ for each $I\in {\cal I}$ and a closed $K\subseteq I$; $K\cap \bar{I%
}\in {\cal K}$ for each $K\in {\cal K}$ and an open $I\subseteq K$. A
topology $({\cal I},{\cal K})$ is {\em Hausdorff}, if the families ${\cal I},%
{\cal K}$ are approximate. If the ring ${\cal B}$ contains a unity $%
E\subseteq X$, then the family ${\cal I}$ of open subsets consists of all
complements (with respect to $E$) $I=E\cap \bar{K}=E\setminus K$ of the
closed subsets $K\in {\cal K}$, and, conversely, the family ${\cal K}$
consists of the complements $K=E\cap \bar{I}=E\setminus I$, $I\in {\cal I}$.
The intersection ${\cal A}={\cal I}\cap {\cal K}\subseteq {\cal B}$ is a
boolean ring of clopen subsets in ${\cal B}$, it contains zero $\emptyset
\in {\cal A}$ and unity $E\in {\cal A}$, if the latter is defined in ${\cal B%
}$.

The notion of topology in a boolean ring ${\cal B}$ is a generalization of
the notion of topological space, which is defined by a dual system ${\cal T}%
=({\cal I},{\cal K})$ of open $I\in {\cal I}$ and closed $K\in {\cal K}$
subsets in the field ${\cal B}={\cal A}(X)$ of all subsets of the set $X$;
here $\cup I_i\in {\cal I}$, $\cap K_i\in {\cal K}$ for arbitrary families $%
\{I_i\}\subseteq {\cal I}$, $\{K_i\}\subseteq {\cal K}$, $X\in {\cal I}\cap 
{\cal K}$ and a set $I$ ($K$) is open (closed) iff the complement $\bar{I}$
is closed ($\bar{K}$ is open): 
\[
{\cal I}=\{\bar{K},\ K\in {\cal K}\},{\cal K}=\{\bar{I},\ I\in {\cal I}\}. 
\]
The smallest $\tau $-field ${\cal A}\subset {\cal A}(X)$, containing all
open and closed subsets of the topological space $(X,{\cal T})$, is called a 
{\em $\tau $-borel field} (a {\em borel field}, if $\tau =\sigma $).

Every topology $({\cal I},{\cal K})$, given on a ring ${\cal B}\subseteq%
{\cal A}(X)$, generates a topological space $(X,{\cal T})$, all unions $\cup
I_i$ and also the whole set $X$ being the open subsets, here the space $(X,%
{\cal T})$ is Hausdorff if the family ${\cal I}$ is attainable and $\sup 
{\cal I}=X$. The $\tau$-borel field ${\cal A}$ contains ${\cal B}$ if the
latter is a $\tau$-field generated by the family $({\cal I},{\cal K})$, e.g.
if ${\cal I}={\cal B}={\cal K}$. More generally, let $f:\ (X)\to X^{\prime}$
be a partial mapping of the topological space $(X,{\cal T})$ into a space $%
X^{\prime}$ with a $\tau$-ring ${\cal B}\subseteq {\cal A}(X^{\prime})$,
generated by a topology $({\cal I},{\cal K})$, and let the preimage $f^*(I)$
of each open subset $I\in {\cal I}$ be open in $(X,{\cal T})$. Then $f$ is a
continuous mapping with respect to the generated topology ${\cal T}^{\prime}$
on $X^{\prime}$, and $h=f^*$ is a $\tau$-homomorphism of the $\tau$-ring $%
{\cal B}$ into the $\tau$-borel field ${\cal A}\subseteq {\cal A}(X)$.

\subsection{Observations and representations}

The observed events are usually described by means of elements of some
semiring (ring) ${\cal B}$. To determine an observation in a system,
described by a logic $\ell$, means to give a logical corresponding of the
observed events $b\in {\cal B}$ by some propositions $a=h(b)\in\ell$. Each
homomorphism $h:\ {\cal B}\to \ell$, describing such representation, is
called an {\em observable} (a {\em quasiobservable}) on a logic (quasilogic) 
$\ell$ over the semiring (ring) ${\cal B}$. The observables on $\tau$-logics
( $\tau$-quasilogics) are described by the $\tau$-homomorphisms of $\tau$%
-semirings ($\tau$-rings) ${\cal B}$, usually, the borel fields - the $%
\sigma $-rings ${\cal B}\subseteq {\cal A}(E)$ of subsets of some set $E$
are studied.

The image $h({\cal B})$ of the semiring ${\cal B}$ generates a boolean ring
(a $\tau $-ring) ${\cal A}\subseteq \ell $ of some propositions in the logic
($\tau $-logic) $\ell $, that commute $a_1\downarrow a_2$ if $a_1,a_2\in 
{\cal A}$. For this reason, it is always possible to describe one fixed
observable on the logic $\ell $ by the representation of the boolean
semiring ${\cal B}$ in some ring ${\cal A}$, usually the boolean algebra $%
{\cal A}\subseteq {\cal A}(X)$ of subsets of some set $X$ is considered.
This choice is based on the Stone theorem on the existence of a
representation of any boolean algebra in some field of subsets, the field of
clopen subsets of some compact, completely disconnected topological space $%
(X,{\cal T})$.

We will show the Stone's construction, generalizing it to the case of the
semirings ${\cal B}$ of finite rank, and also of the boolean ring, not
necessarily containing unity. Let us denote $X_0$ the set of the maximal
filters of the semiring ${\cal B}$, note that the filter $x=\nabla $,
defined, in the boolean semiring, by the conditions: 
\[
1)\ 0\notin \nabla ,\quad 2)\ a\in \nabla ,\ b\geq a\Rightarrow b\in \nabla
,\quad 3)\ a,b\in \nabla \Rightarrow a\wedge b\in \nabla 
\]
is maximal iff for any $a\notin \nabla $ there is a $b\in \nabla $, such
that $a\wedge b=0$. Assigning the set $h_0(b)=\{x\ :\ b\in x\}$ to each $%
b\in {\cal B}$, we get by this definition: 
\[
h_0(0)=\emptyset ,\quad b\geq a\Rightarrow h_0(b)\supseteq h_0(a),\quad
h_0(a\wedge b)=h_0(a)\cap h_0(b). 
\]
Moreover, for any finite family $\{a_k\}$, from the condition $a_k\notin x$
for each $k$ follows the existence of 
\[
b=\wedge b_k,b_k\in x,b_k\wedge a_k=0, 
\]
for which $b\wedge \left( \vee a_k\right) =0$, by distributivity, i.e. $\vee
a_k\notin x$, or $x\notin h_0(\vee a_k)$, if $x\notin h_0(a_k)$ for each $k$.

Because, on the other hand, from $a_k\leq \vee a_k$ follows that $%
h_0(a_k)\subseteq h_0(\vee a_k)$, we have $h_0(\vee a_k)=\cup h_0(a_k)$ for
any finite family $\{a_k\}$. Hence the mapping $h_0:{\cal B}\to {\cal A}%
(X_0) $ is a homomorphism of the semiring ${\cal B}$ into the field of
subsets of maximal filters $x\in X_0$, preserving the disjointness: 
\[
a\wedge b=0\Rightarrow h_0(a)\cap h_0(b)=\emptyset 
\]
and, consequently, the operation of subtraction 
\[
h_0(b\setminus a)=h_0(b)\cap \bar{h_0(a)}, 
\]
if ${\cal B}$ is a ring. Such representation is faithful, because for each $%
a\neq 0$, there is at least one maximal filter $x$: $a\in x$ (the filter $%
\nabla =\{b\ :\ b\geq a\}$ can be extended to a maximal filter, by the axiom
of choice). The semiring (ring) ${\cal B}_0=h_0({\cal B})$ has a finite rank
and is perfect, because each maximal filter is $\nabla \subset {\cal B}_0$
is a filter $\{h_0(b)\ :\ b\in x\}$, determined by some element $x\in X_0$, 
\[
\nabla =\{B\in {\cal B}_0\ :\ x\in B\}, 
\]
and separating: for each $x_1\neq x_2$ there is a $b\in x_1$, $b\notin x_2$,
i.e. $x_1=h_0(b)$, $x_2\notin h_0(b)$. Moreover, every $B\in {\cal B}_0$ is
an open compact subset of the locally compact completely disconnected space $%
(X_0,{\cal T}_0)$ with the base ${\cal B}_0$, and each open compact subset $%
A\subseteq X_0$ in the topology ${\cal T}_0$ is a finite sum $%
A=\sum_{k=1}^nB_k$ of disjoint $B_k\in {\cal B}_0$.

First, each $B\in {\cal B}_0$ is open, because according to the definition
of the topology ${\cal T}_0$, the open subsets $A$ of $X_0$ are all the
unions $A=\cup B_i$, $B_i\in {\cal B}_0$. On the other hand, each $B\in 
{\cal B}_0$ is closed, as $\bar{B}$ is open: 
\[
\bar{B}=\bar{B}\cap X_0=\cup (\bar{B}\cap B_i)=\cup B_{i_k}, 
\]
where $B_i$ cover $X_0=\cup B_i$ and $B_{i_k}\in {\cal B}_0$ is a finite
decomposition of $\bar{B}\cap B_i=\cup _{k=1}^{n_i}B_{i_k}$ (the set of
finite unions $\cup B_k$, $B_k\in {\cal B}_0$ is a boolean ring, invariant
under the operation $A\setminus B=\bar{B}\cap A$). Consequently, each $B\in 
{\cal B}$ is clopen and compact, because it is a space of Stone's
representation of the boolean algebra ${\cal A}=\{a\in {\cal B},\ a\leq b\}$%
, $h_0(b)=B$. As the sets $B\in {\cal B}_0$ cover $X_0$ and separate its
points, the space $(X_0,{\cal T}_0)$ is locally compact and completely
disconnected, and every open compact set $A\subset X_0$ of this topological
space has the form $A=\sum_{k=1}^nB_k$ (the covering $A=\cup B_i$ has a
finite subcovering, that can be decomposed into a sum of disjoint $B_k\in 
{\cal B}_0$ by the finiteness of the rank of the semiring ${\cal B}$). Note
that if ${\cal B}$ is a boolean ring, then the representing semiring ${\cal B%
}_0=h_0({\cal B})$ is identical with the generated ring ${\cal A}_0$ of all
open compact subsets of the space $(X_0,{\cal T})$. As the space $X_0=h_0(1)$
is compact if $1\in {\cal B}$, the ring ${\cal A}_0=h_0({\cal A})$ for any
boolean algebra ${\cal A}$ is a boolean algebra of clopen subsets with the
complementation: 
\[
\overline{h(a)}=h(1-a)=h(\bar{a}) 
\]

Unfortunately, the isomorphism $h_0:{\cal B}\rightarrow {\cal B}_0$ is only
finite-additive and cannot be used for the construction of $\tau $%
-representations in the category of semirings of the infinite rank $\tau $ .
The maximal $\tau $-subrepresentation $h_\tau :{\cal B}\rightarrow h_0\left(
b\right) \cap X_\tau $ of the universal representation $h_0$ , which
represents the elements $b\in {\cal B}$ by the subsets $B\subseteq X_\tau $
of maximal $\tau$-filters $x\in X_\tau $ , is not in general faithful. It is
faithful only if for each $b\neq 0$ there exists a maximal $\tau $-filter $%
\nabla $ containing $b$ , what cannot be true already for $\tau =\sigma $ .
In general case the universal $\tau $-representation can be obtained as a
factor-representation $j_\tau \circ h_0:{\cal B}\rightarrow {\cal A}_\tau
/\Delta _\tau $ in the $\tau $-completion ${\cal A}_\tau ={\cal A}_0\vee
\Delta _\tau $ of the ring ${\cal A}_0$ of the space $X_0$ , factorised with
respect to the $\tau $-ideal of all compact subsets $A\in {\cal A}_\tau
/\Delta _\tau $.of the $\tau $-category. If the ideal $\Delta _\tau $ is
regular, the factor-representation is equivalent to the $\tau $%
-subrepresentation $h_\tau :{\cal B}\rightarrow {\cal A}\left( X_\tau
\right) $ . Although in general case the$\tau $-ideal $\Delta _\tau $ is
smaller than the regular ideal $\left\{ A\in {\cal A}_\tau :A\cap X_\tau
=\emptyset \right\} $ , the universal $\tau $-homomorphism $j_\tau \circ h_0$
could be faithful even in this case as it follows from the next theorem for $%
\tau =\sigma $ .

\begin{theorem}
Each $\sigma $-semiring ${\cal B}$ is $\sigma $-isomorphic with some
semiring ${\cal B}_\sigma \subseteq {\cal A}_\sigma /\Delta _\sigma $ of a
factor algebra of the $\sigma $-completion ${\cal A}_\sigma ={\cal A}_0\vee
\Delta _\sigma $ of the ring all open compact subsets of a locally compact
space $X_0$ of maximal filters $x\subset {\cal B}$ by the $\sigma $-ideal $%
\Delta _\sigma $ of the compact subsets $A\in {\cal A}$ of the first
category. \label{th:2.3}
\end{theorem}

The proof of this theorem is a modification of the proof, given in [4] for
the case of $\sigma $-algebras, to the case of the semirings (rings), which
is analogical to the modification of the Stone theorem, given above.

\subsection{Regular distributions}

If $h:\ {\cal B}\to \ell$ is an observable (quasiobservable) on a logic
(quasilogic) $\ell$, and $m:\ \ell\to R_+$ is a measure on $\ell$, then the
composition $n=m\circ h$ is a distribution on ${\cal B}$, which is $\tau$%
-additive if $m$ and $h$ are. In particular, if ${\cal B}$ is a semiring of
events with unity $1\in {\cal B}$ and a state $p:\ell\to [0,1]$ is given,
then $w=p\circ h$ is a probability distribution on ${\cal B}$, if $h(1)$
belongs to the support of this state, i.e. $w(1)=p(h(1))=1$. Such
distribution, describing the probabilities of observable events $b\in {\cal B%
}$, defines a {\em statistical state} on the semiring ${\cal B}$,
corresponding with the measurement $h$ in a system, described by the logic $%
\ell$ with the given state $p$.

In general, an arbitrary $\tau $-distribution $n:\ {\cal B}\to R_{+}$ on a
semiring ${\cal B}$ is called $\tau $-{\em representable}, if there exist
some $\tau $-homomorphism $h:\ {\cal B}\to {\cal A}$ into some $\tau $-ring $%
{\cal A}\subseteq {\cal A}(X)$ and a $\tau $-measure $\mu :\ {\cal A}\to
R_{+}$, such that $n=\mu \circ h$. According to the theorem on the
representation of boolean semirings (the Stone theorem for boolean
algebras), proved in the preceding section, each distribution $n$ is
finitely-representable on some locally compact space $X$, that can be chosen
as the space $X_0$ of the universal representation $h_0:\ {\cal B}\to {\cal A%
}_0$, extending the finite-additive distribution $\nu _0=n\circ h_0^{-1}$
from the semiring ${\cal B}_0=h_0({\cal B})$ to the measure $\mu _0:\ {\cal A%
}_0\to R_{+}$ on the ring ${\cal A}_0$, by additivity. However, by $\tau
=\sigma $, not every $\tau $-distribution is $\sigma $-representable in the
described sense, because, in general, no universal $\sigma $-isomorphism of
an arbitrary $\sigma $-ring ${\cal B}$ into some field ${\cal A}(X)$ exists.
It is therefore natural to generalize the notion of representability of a $%
\tau $-distribution $n:{\cal B}\to R_{+}$, calling it {\em %
factor-representable}, if a $\tau $-measure $\mu :\ {\cal A}\to R_{+}$ and a 
$\tau $-homomorphism $h:\ {\cal B}\to {\cal A}/\Delta $ into a
factor-algebra of some $\tau $-ring ${\cal A}\subseteq {\cal A}(X)$ by some $%
\tau $-ideal $\Delta \subset {\cal A}$ of subsets $A\in \Delta $ of zero
measure $\mu (A)=0$ exist, such that $n=m\circ h$, where $m$ is a $\tau $%
-measure on ${\cal A}/\Delta $, determining $\mu =m\circ j$ by the
factorisation homomorphism $j:\ {\cal A}\to {\cal A}/\Delta $. According to
the theorem~\ref{th:2.3}, any $\sigma $-distribution is factor-representable
on some locally compact space $X$ with a $\sigma $-measure $\mu $, that
vanishes on a $\sigma $-ideal of subsets of the first category, because as $%
(X,\mu )$ it is always possible to choose the space $X_0$ of the universal
factor-representation $h_\sigma :\ {\cal B}\to {\cal A}_\sigma /\Delta
_\sigma $ and the $\sigma $-measure $\mu _\sigma :\ {\cal A}_\sigma \to
R_{+} $, that is the extension of the measure $\mu _0$ from the ring ${\cal A%
}_0$ to the $\sigma $-complete ring ${\cal A}_\sigma ={\cal A}_0\vee \Delta
_\sigma $.

For $\tau>\sigma$, $\tau$-distributions $n:\ {\cal B}\to R_+$ exist that are
not factor-representable, even if ${\cal B}$ is a boolean $\tau$-algebra;
however, if $n$ is factor-representable, then the space $X$ can be always
chosen locally compact, for example, $X_0$ with the smallest $\tau$-ring $%
{\cal A}_{\tau}\supseteq{\cal A}_0$ and the $\tau$-ideal $%
\Delta_{\tau}\subseteq {\cal A}_{\tau}$, putting $h=j_{\tau}\circ h_0$. The
corresponding $\tau$-measure $\mu_{\tau}:\ {\cal A}_{\tau}\to R_+$ is the
extension by $\tau$-additivity of the measure $\mu_0$ from the ring ${\cal A}%
_0$ to the $\tau$-complete ring ${\cal A}_{\tau}={\cal A}_0\vee\Delta_{\tau}$%
, that is, as in the case $\tau=\sigma$, completely regular $\mu_{\tau}=\sup
\{ \mu_{\tau}(K)\ :\ K\subseteq A,\ k\in {\cal K}\}$ for any $A\in {\cal A}$%
, with respect to the family ${\cal K}$ of the compact subsets of the space $%
X_0$. The last follows from the condition of $\tau$-approximity of the
family ${\cal K}_{\tau}\subseteq{\cal K}$ of the $\tau$-closed compact
subsets $K\in {\cal K}\cap{\cal A}_{\tau}={\cal K}_{\tau}$ of the ring $%
{\cal A}_{\tau}$.

In general, if $X$ is a locally compact space and $\mu :\ {\cal A}\to R_{+}$
is a measure on some boolean ring ${\cal A}\subset {\cal A}(X)$ with the
property of approximity with respect to some family ${\cal K}$ of compact
subsets $K\in {\cal A}$, then $\mu $ is a $\sigma $-measure and it can be
always extended to a regular measure $\bar{\mu}$ on the completion $\bar{%
{\cal A}}={\cal A}\vee \bar{\Delta}$ of the ring ${\cal A}$ with the compact
subsets $K\in \bar{\Delta}$ of measure $0:$ 
\[
\bar{\mu}(K)=\inf \{\mu (A)\ :\ A\supset K,\ A\in {\cal A}\}. 
\]

\section{Quantum systems}

\subsection{Hilbertian clans}

Quantum systems, which are, strictly speaking, any physical systems, due to
the uncertainty principle, are usually described by a Hilbert space ${\cal E}
$, the subspaces of which ${\cal H}\subseteq {\cal E}$ form a {\em quantum
logic} with respect to inclusion: ${\cal H}_1\leq {\cal H}_2$ if ${\cal H}%
_1\subseteq {\cal H}_2$, and orthogonal complement $\bar{{\cal H}}={\cal H}%
^{\perp }$. The quantum logic ${\cal L}({\cal E})$ of all subspaces ${\cal H}%
\in {\cal L}({\cal E})$ has a zero - the one-dimensional (zero-dimensional)
subspace ${\cal O}=\{0\}$, a unit - the whole space ${\cal E}\in {\cal L}(%
{\cal E})$, it contains the greatest lower bound 
\[
\bigwedge_{i\in I}{\cal H}_i=\bigcap_{i\in I}{\cal H}_i, 
\]
and, consequently, also the least upper bound 
\[
\bigvee {\cal H}_i=(\cap {\cal H}_i^{\perp })^{\perp } 
\]
for an arbitrary family ${\cal H}_i\subseteq {\cal E},i\in I$. The latter
logical structure is nondistributive, except in the trivial case of the
one-dimensional ${\cal E}$, as can be easily seen taking into account any
nonorthogonal subspaces ${\cal H}_1,{\cal H}_1$ intersecting only in zero, $%
{\cal H}_1\wedge {\cal H}_2=0$, which describe nonsimultaneous but not
disjoint propositions, for example, eigensubspaces ${\cal H}_{\Delta p},%
{\cal H}_{\Delta q}$, of the momentum operator $\hat{p}$ and the coordinate
operator $\hat{q}$ of a quantum mechanical particle, corresponding to finite
values $p\in \Delta p$, $q\in \Delta q$. In particular, this indicates that
the quantum logic ${\cal L}({\cal E})$ is not equivalent, by dim${\cal H}>1$%
, to any boolean algebra, and it cannot be represented by any field of
subsets of a set $X$, i.e., a quantum system cannot be described in a
classical way in terms of a phase space $X$.

From the algebraic point of view, it is appropriate to identify propositions
of a quantum logic ${\cal L}({\cal E})$ not with the subspaces ${\cal H}%
\subseteq {\cal E}$, but with the corresponding orthogonal projections $%
P:\xi \in {\cal E}\mapsto P\xi \in {\cal H}$ onto the subspace ${\cal H}=P%
{\cal E}$, which are characterized, in the algebra ${\cal B}({\cal E})$ of
all bounded linear operators $A$ in ${\cal E}$, as selfadjoint idempotents 
\[
A^{*}=A=A^2. 
\]
The logical partial operations of substraction, addition and multiplication
of propositions are described by the algebraic operations $A-B,A+B,AB$,
defined in the logic ${\cal L}({\cal E})$ of all orthoprojections of the
space ${\cal E}$ for comparable in the sense $AB=B$, orthogonal in the sense 
$AB=0$, and commuting in the sense $AB=BA$, couples $A,B\in {\cal L}({\cal E}%
)$, respectively. We note that also any boolean algebra can be represented
in this way, if we assigne to the elements $a\in {\cal A}$ e.g. the
orthoprojectors on the Hilbert subspaces ${\cal H}_a\subset {\cal E}$ of
square summable functions $x\mapsto \xi (x)$ on the set $X$ of maximal
ideals of the algebra ${\cal A}$ with carriers $h(a)\subseteq X$, defined by
the Stone representation $a\mapsto h(a)$.The corresponding subset of
orthoprojections $P_a\in {\cal L}({\cal E})$ is commutative, and forms a
boolean ring with unit with respect the operations of conjunction $A\wedge
B=AB$, difference $A\setminus B=A-AB$ and disjunction 
\[
A\vee B=AB+A\setminus B+B\setminus A. 
\]
It is not difficult to see that this ring may be defined as the set of
selfadjoint idempotents of a commutative subalgebra ${\cal A}\subset {\cal B}%
({\cal H})$ of operators of pointwise multiplication $A\xi (x)=\alpha (x)\xi
(x)$ by bounded number valued functions $x\mapsto \alpha (x)$.

In general, we will call a clan an arbitrary logic ${\cal L}\subset {\cal L}(%
{\cal E})$ of orthoprojections, consisting of selfadjoint idempotents of any 
$*$-algebra ${\cal A}\subset {\cal B}({\cal E})$ of bounded operators $A\in 
{\cal A}\Rightarrow A^{*}\in {\cal A}$ containing an orthoprojector $P\neq 0$
as the algebraic unit: 
\[
PA=A=AP,\forall A\in {\cal A}. 
\]
A clan ${\cal L}$ is $\tau $-complete if ${\cal A}$ is a $\tau $-complete
algebra, that is, if any downward directed set $\{A_i\}\subset {\cal A}$ of
positive operators $A_i\geq 0$ has in ${\cal A}$ the infimum $\bigwedge
A_i\in {\cal A}$, and is complete if ${\cal A}$ is a complete $*$-algebra,
for example, a von Neumann algebra. A Banach $\sigma $-complete $*$-algebra $%
{\cal A}$ is called a Baer $*$-algebra; the corresponding clan ${\cal L}%
\subset {\cal A}$ consists of carriers of the operators $A\in {\cal A}$,
i.e., the smallest orthoprojections in ${\cal E}$, $E\in {\cal L}({\cal E})$%
, such that $EA=A=AE$.

The following theorem shows that if any pair $P_1, P_2$ of orthoprojections
in a clan ${\cal L}$ for which $P_1\wedge P_2 = 0$ is orthogonal, then the
clan is a boolean algebra.

\begin{theorem}
A necessary and sufficient condition of the distributivity of a clan ${\cal L%
}$ is the equivalence of the conditions of nonsimultaneousness $P_1\wedge
P_2=0$ and mutual exclusivness $P_1P_2=0$. \label{th:2.4}
\end{theorem}

\noindent{\sc Proof} Evidently, $P_1P_2=0$ implies $P_1\wedge P_2=0$, so
that it suffices to prove the implication 
\[
P_1\wedge P_2=0\Rightarrow P_1P_2=0, 
\]
i.e., $P-P_1\geq P_2$ whenever $P\geq P_1,P_2$, from distributivity. We have 
\[
P_2=P\wedge P_2=((P-P_1)\vee P_1)\wedge P_2=((P-P_1)\wedge P_2)\vee
(P_1\wedge P_2)=(P-P_1)\wedge P_2. 
\]
Conversely, writing any pair $P_1,P_2\in {\cal L}$ in the form $%
P_i=P_1\wedge P_2+B_i$, where $B_i=P_i-P_1\wedge P_2$, $i=1,2$, we get 
\[
P_1P_2=P_1\wedge P_2=P_2P_1, 
\]
if the implication $P_1\wedge P_2=0\Rightarrow P_1P_2=0$ is true. From this
the distributivity follows: 
\begin{eqnarray*}
&&(P_1\vee P_2)\wedge P_3=(P_1P_2+B_1+B_2)P_3 \\
&=&P_1P_2P_3+B_1P_3+B_2P_3=P_1P_3\vee P_2P_3.
\end{eqnarray*}

From the proof of Theorem~\ref{th:2.4} it also follows that if the
suppositions of the theorem are satisfied, the set ${\cal L}$ of
orthoprojections is commutative and, therefore, it is a clan of some
commutative algebra ${\cal A}$.

\subsection{Quantum states}

If ${\cal L}({\cal E})$ is a quantum logic of orthoprojections in a Hilbert
space ${\cal E}$ with a scalar product $(\xi ,\xi )$, then for any vector $%
\xi \in {\cal E}$, the formula $\rho ^\xi (A)=(A\xi ,\xi )$ defines a
positive completely additive mapping $\rho ^\xi :{\cal L}\to [0,1]$, which
is a state on ${\cal L}$, if $P\xi =\xi $ for the unit $P$ of the logic $%
{\cal L}$. In general, any linear functional $\rho $ on a $*$-algebra ${\cal %
A}\subset {\cal B}({\cal E})$ with a unit $P$, satisfying condition of
positivity $\rho (A)\geq 0$ whenever $a\geq 0$, is called a distribution on $%
{\cal A}$, since it defines a measure $m(A)=\rho (A)$ on the clan ${\cal L}$
of orthoprojections of the algebra. The measure $m$ is bounded: $\Vert
m\Vert =m(P)$, and it is a $\tau $-measure if the functional $\rho $
satisfies the condition $\wedge _{i\in I}\rho (A_i)=0$ for any decreasing $%
A_i\downarrow 0$ family $\{A_i\}_{i\in I}$, $|I|\leq \tau $, of operators $%
A_i\in {\cal A}$. A distribution $\rho $ is called a vector distribution if $%
\rho (A)=(A\xi ,\xi )=\rho ^\xi (A)$ for some vector $\xi \in {\cal E}$. Any
convex combination of vector distributions is a normal distribution, which
defines a normal measure on ${\cal L}$, because then $\bigwedge_{i\in I}\rho
(A_i)=0$ if $A_i\downarrow 0$ for any $I$.

We will consider now an arbitrary (abstract) ``algebra of observables'' $%
{\cal A}$, defined as complex linear space with an operation of associative
multiplication 
\[
(a,b)\mapsto a.b,a(b+c)=ab+ac,a(\lambda b)=\lambda ab,\forall \lambda \in 
{\bf C}, 
\]
and an involution $a\mapsto a^{*}$, $a^{**}=a$, for which 
\[
(a+b)^{*}=a^{*}+b^{*},(ab)^{*}=b^{*}a^{*},(\lambda a)^{*}=\bar{\lambda}%
a^{*}. 
\]
Assume that in ${\cal A}$, a relation of partial ordering is defined,
satisfying the conditions 
\[
a\leq c\Rightarrow a+b\leq c+b,b^{*}ab\leq b^{*}cb;\lambda a\leq \lambda c 
\]
for any $b\in {\cal A}$ and $\lambda \geq 0$, and let for any $a\in {\cal A}$
there is a self adjoint idempotent $e\geq 0$ in ${\cal A}$ such that $%
ea=a=ae $. Since $a^{*}a=a^{*}ea\geq 0$ for any $a\in {\cal A}$, the
positive cone ${\cal A}_{+}=\{a\geq 0\}$ of the algebra ${\cal A}$ contains
all elements of the form $a^{*}a$, in particular, all selfadjoint
idempotents, which form an upward directed subset $\ell $ of the cone ${\cal %
A}_{+}$, since for any $e_1,e_2\in \ell $ there exists $e\in \ell $, $e\geq
e_i$, $i=1,2$, determined by the condition $ea=a=ae$ for $a=e_1+\mbox{i}e_2$%
. It is not difficult to verify that the set $\ell \subset {\cal A}$ is a
logic with the order relation $e_1\leq e_2$ if $e_1e_2=e_1$, commutation
relation $e_1\bar{\downarrow}e_2$ if $e_1e_2=e_2e_1$, and orthogonality
relation $e_1\perp e_2$ if $e_1e_2=0$, without a logical unit, in general,
which is $\tau $-complete, if the cone ${\cal A}_{+}$ is $\tau $-complete in
the sense $\bigwedge_{i\in I}a_i\in {\cal A}_{+}$ for any decreasing set $%
\{a_i\}_{i\in I}\subset {\cal A}_{+}$ of cardinality $|I|\leq \tau $.

A distribution on an algebra of observables ${\cal A}$ is any linear
functional $r:{\cal A}\to {\bf C}$, satisfying condition of positivity $r(a)
\geq 0$ whenever $a\geq 0$; if in addition, $\wedge r(a_i) = 0$ for any $%
\tau $-family $\{ a_i\} \subset {\cal A}_+$ decreasing to zero, then $r$ is
called a $\tau$-distribution (normal distribution if the latter condition
holds for any cardinal $\tau$).

Distributions on an algebra ${\cal A}$, with a logic $\ell \subset {\cal A}$
as an approximative unit, are defined by their values on the set of elements
of the form $a^{*}a$, since every self adjoint element $a=a^{*}$ can be
written in the form $a=a_{+}-a_{-}$, where $a_{\pm }=\frac 14(a\pm e)^2$,
where $e^{*}=e=e^2\in \ell $ satisfies the condition $ea=a=ae$. The
restriction $m=r|\ell $ defines, evidently, a measure on $\ell $, which is a 
$\tau $-measure if $r$ is a $\tau $-distribution. A state on ${\cal A}$ is
any distribution $r$, iducing a state on $\ell $, i.e., a probability
measure, for which $r(e)=1$ for some $e\in \ell $. If all elements $a=a^{*}$
are normed by the functional 
\[
\Vert a\Vert =\mbox{inf}\{r>0:re>\pm a\}, 
\]
where $e\in \ell $ is any idempotent, satisfying the condition $eae=a$, then
by Hahn-Banach theorem there exists a sufficiently rich set of states to
satisfy the condition $r(a)=1$ for any $a^{*}=a$ for which $\Vert a\Vert =1$%
, and, in particular, the condition $r(e)=1$ for any nonzero $e\in \ell $.
This condition is satisfied by any normed algebra of observables ${\cal A}$
such that $\Vert a\Vert =\Vert aa^{*}\Vert ^{1/2}$, called a C$*$-algebra if
it is a Banach algebra, in particular, any $*$-subalgebra ${\cal A}\subset 
{\cal B}({\cal E})$ of bounded operators, acting on some pre-Hilbert space $%
{\cal E}$. As it follows from the construction in the following paragraph,
this case actually describes all algebras of observables ${\cal A}$ normed
by elements $e\in \ell $.

\subsection{Quantum representations}

A representation of a partially ordered $*$-algebra ${\cal A}$ is any $*$%
-homomorphism of this algebra into an operator $*$-algebra ${\cal B}({\cal E}%
)$ of some pre-Hilbert space ${\cal E}$, that is, a linear mapping $\pi:%
{\cal A}\to {\cal B}({\cal E})$, for which 
\[
\pi(ab) = \pi(a)\pi(b),\ \pi(a^*) = \pi(a)^*, \ a\geq 0 \Rightarrow
\pi(a)\geq 0. 
\]
Any representation $\pi$ defines a homomorphism of the logic $\ell\subset 
{\cal A}$ into the logic ${\cal L}({\cal E})$ of all orthoprojections of the
space ${\cal E}$, which is a $\tau$-homomorphism (normal homomorphism) if $%
\bigwedge_{i\in I} \pi(a_i) = 0$ for any decreasing to 0, $a_i\downarrow 0$,
family $\{a_i\}$ of cardinality $\tau$ (of any cardinality). It turns out
that any logical state on $\ell$, defined by an algebraic state $r:{\cal A}%
\to {\bf C}$, is described in the corresponding representation by a normed
vector $\xi_1\in {\cal E}$, and this representation is a $\tau$%
-representation (a normal representation) if $r$ is a $\tau$-state (a normal
state). We will bring a construction of such representation, which is a
generalization of the well-known GNS-construction to the case of an
arbitrary $*$-algebra ${\cal A}$.

Let $r$ be a state on ${\cal A}$ and element $e_1\in \ell $ satisfy
condition $r(e_1)=1$. Denote by ${\cal E}$ the set of classes 
\[
\xi =\{a^{\prime }:a^{\prime }\simeq a^{*}\}\equiv |a\rangle 
\]
of elements $a^{\prime }$ equivalent to $a^{*}$ with respect to the kernel $%
{\cal N}=\{a:r(aa^{*})=0\}$ of the hermitian form $\langle b|a\rangle
=r(ba^{*})$ on ${\cal A}$, where we put $a^{\prime }\simeq a^{*}$ if $%
a^{\prime }-a^{*}\in {\cal N}$. This set forms a pre-Hilbert space with
respect to the scalar product $(\xi ,\eta )=\langle b|a\rangle =\eta ^{*}\xi 
$, where $\eta ^{*}\equiv \langle b|$ is a linear functional $|a\rangle
\mapsto r(ba^{*})$, corresponding to the vector $\eta =|b\rangle $. To every 
$b\in {\cal A}$, we define a linear operator $\pi _r(b):|a\rangle \mapsto
|ab^{*}\rangle $, having an adjoint $\pi _r(b)^{*}=\pi _r(b^{*})$ on ${\cal E%
}$, in correspondence with the equality 
\begin{eqnarray*}
(\pi _r(b)\xi ,\xi ) &=&\langle a|ab^{*}\rangle =r(aba^{*}) \\
&=&\langle ab|a\rangle =(\xi ,\pi _r(b^{*})\xi ),\ \forall \xi \in {\cal E}.
\end{eqnarray*}

The mapping $\pi _r:b\to \pi (b)$ is, evidently, linear, multiplicative 
\[
\pi _r(b_1b_2)|a\rangle =|a(b_1b_2)^{*}\rangle =|ab_2^{*}b_1^{*}\rangle =\pi
_r(b_1)\pi _r(b_2)|a\rangle , 
\]
preserves order, and $\wedge (\pi _r(a_i)\xi ,\xi )=0$ whenever $%
a_i\downarrow 0$, $|I|\leq \tau $, and $r$ is $\tau $-normal. In this way,
we construct a $\tau $-normal representation of an algebra ${\cal A}$ with
an approximative unit $\ell \subset {\cal A}$. It remains to find a vector
in ${\cal E}$, describing the given state on $\ell $. We prove that we can
take any vector $\xi _1=|e_1\rangle $, where $e_1\in {\cal A}$ is a
self-adjoint idempotent, satisfying the condition $r(e_1)=\mbox{sup}_{e\in
\ell }r(e)$,i.e., we have the representation 
\[
r(a)=\langle e_1|\pi _r(a)|e_1\rangle =r(a_1ae_1) 
\]
for any $a\in {\cal A}$.

Indeed, due to the directedness of the set $\ell $, for any $a\in {\cal A}$
there exists a self-adjoint idempotent $e\in \ell $, majorating $e_1$, such
that $ea=a=ae$. It follows that 
\[
r(a)=r(e_0ae_0)+r(e_0ae_1)+r(e_1ae_0)+r(e_1ae_1), 
\]
where $e_0=e-e_1\in \ell $ is a self-adjoint idempotent satisfying $r(e_0)=0$
owing to $r(e_1)=r(e)$. But this means that 
\[
r(e_0ae_0)=r(e_0ae_1)=r(e_1ae_0)=0, 
\]
due to the generalized Schwartz inequality 
\[
|r(ba^{*})|^2\leq r(bb^{*})r(aa^{*}), 
\]
which can be derived from the Sylvester determinant criterion of
non-negativeness of the 2x2-form 
\[
r((\alpha a+\beta b)(\alpha a+\beta b)^{*})\geq 0 
\]
for any $\alpha ,\beta \in {\bf C}$ applied to $b=e_0$. In this way, we
proved the following

\begin{theorem}
Let $\ell $ be a logic of self-adjoint idempotents, which is an
approximative unit of a $*$-algebra ${\cal A}$, and $\rho :\ell \to [0,1]$
be a probability state on $\ell $, defined by some $\tau $-normal
distribution $r$ on ${\cal A}$, $\rho =r/\ell $. Then there exist a
pre-Hilbert space ${\cal E}$, a $\tau $-normal representation $\pi :{\cal A}%
\to {\cal B}({\cal E})$ and a normed vector $\xi _1\in {\cal E}$, which
defines a state as a vector state on the clan of the orthoprojectors 
\[
E=\pi (e):\rho (e)=(\pi (e)\xi _1,\xi _1). 
\]
\label{th:2.5}
\end{theorem}

Let us note that in case that an injective representation $\pi$ exists,
obtained, e.g., by the described construction for a strictly positive
(faithful) distribution $r$, the quasilogic $\ell$ of the algebra ${\cal A}$
is isomorphic to the clan ${\cal L}$ of the operator algebra {\sl A} $= \pi(%
{\cal A})$, and can be identified with it. In particular, for any algebra of
observables ${\cal A}$, normed by the logic $\ell \subset {\cal A}$, there
exists a faithful $\tau$-representation $\oplus \pi_{r}$ in the Hilbert
space $\oplus {\cal E}_{r} = {\cal H}$, the completion of the direct sum $%
\oplus {\cal E}_{r}$ over all states $r: {\cal A}\to {\bf C}$.

\subsection{Quantum observables}

Let ${\cal L}$ be a clan of orthoprojections, representing a quantum logic
of a system in a Hilbert space ${\cal E}$ of vectors $\xi \in {\cal E}$ of
quantum states, and ${\cal H}$ be a Hilbert space, describing another
quantum system, which is accessible to observations. Any mapping $\rho ^F:%
{\cal L}\to {\cal B}({\cal H})$ of the form 
\[
\rho ^F(A)=F^{*}AF, 
\]
where $F:{\cal H}\to {\cal E}$ is a linear operator whose adjoint on ${\cal E%
}$ is $F^{*}$, $F^{*}{\cal E}\subseteq {\cal H}$, will be called an operator
distribution on the clan ${\cal L}$ with respect to the Hilbert space ${\cal %
H}$ . A distribution $\rho $ is called an operator state if $F$ is an
isometric embedding of ${\cal H}$ into ${\cal E}$: $F^{*}F=I$ and $%
FF^{*}\subseteq P$ for some orthoprojector $P\in {\cal L}$. In the case of
one-dimensional space ${\cal H}={\bf C}$, the operator distribution (state)
coincides with the vector distribution (state) defined by the
one-dimensional operator $F=\xi :\alpha \in {\bf C}\mapsto \alpha \xi \in 
{\cal E}$. Any operator distribution $\rho ^F$ on ${\cal L}$ defines a
positive completely additive normal mapping 
\[
A\in {\cal L}\mapsto \rho ^F(A)\geq 0, 
\]
which maps some orthoprojection $P\in {\cal L}$ to the identity operator $I$
of ${\cal H}$: 
\[
\rho ^F(P)=F^{*}PF=I. 
\]
A triple $({\cal E},{\cal L},F)$, where $F$ is an operator of a state $\rho
^F$ on ${\cal L}$, will be called a quantum operator- probability space
(quantum probability system), and simply a quantum probability space if $%
F=\xi $ is a vector of the state (${\cal H}$ is one-dimensional).

Let ${\cal B}$ be a semiring (ring) with unit $E$, describing events $B\in 
{\cal B}$ that can be observed in a quantum system $({\cal E},{\cal L},F)$
``on the scale of the measuring device $E$''. Any homomorphism $h:{\cal B}%
\to {\cal L}$, for which $h(E)F=F$, induces for any normed vector $\eta \in 
{\cal H}$, a probability state on ${\cal B}$ described by the measure 
\[
m^\eta (B)=(h(B)F\eta ,F\eta ), 
\]
and is called a quantum observable in the system $({\cal E},{\cal L},F)$.
The mapping 
\[
B\mapsto m^F(B)=F^{*}h(B)F, 
\]
defined by a quantum observable $h$, evidently, is positive and additive,
satisfies condition $m^F(B)\leq I$, where $I$, where $I$ is the unit
operator in ${\cal H}$. The homomorphism from a Boolean ring ${\cal B}$ into
the quasilogic $[0,I]\subset {\cal B}({\cal H})$ of positive operators $%
A\leq I$ on the Hilbert space ${\cal H}$, satisfying the normalization
condition $m^F(E)=I$, is called {\em quasiobservable} over ${\cal B}$. Any
quasiobservable $m:{\cal B}\to [0,I]$, which defines probability
distributions $m^\eta (B)=(M(B)\eta ,\eta )$ corresponding to unit vectors $%
\eta \in {\cal H}$, will be called an operator-valued probability
distribution (o.p.d.) on ${\cal B}$, or a weak quantum observable in the
Hilbert space ${\cal H}$. Quantum observables in the category of $\tau $%
-semirings ${\cal B}$ are described by $\tau $-homomorphisms $h$ into the
operator $\tau $- clans ${\cal L}\subset {\cal B}({\cal E})$, inducing
operator-valued probability $\tau $-distributions 
\[
m^F(B)=F^{*}h(B)F 
\]
as the weak quantum observables in ${\cal H}$, and they are normal in case
of normal $h$.

Two observables $h_1, h_2$ in possibly different quantum systems $({\cal E}%
_i,{\cal L}_i,F_i)$, $i=1,2$ are said to be equivalent with respect to the
Hilbert space ${\cal H}$ (or ${\cal H}$-equivalent) if they induce the same
weak quantum observable, i.e., if $m^{F_1} = m^{F_2}$. The existence and
canonical construction of the class of equivalent observables, corresponding
to any o.p.d. $m$, is stated by the following theorem, which generalizes the
Naimark construction in the case of arbitrary (abstract) $\tau$-semirings $%
{\cal B}$ with a unit $E\in {\cal B}$.

\begin{theorem}
Let $m:{\cal B}\to {\cal B}({\cal H})$ be an operator-valued distribution on
a $\tau $-semiring ${\cal B}$ with a unit $E\in {\cal B}$ satisfying the
following conditions:

\begin{enumerate}
\item  positivity: $(m(B)\eta ,\eta )\geq 0$, $\forall B\in {\cal B}$, $\eta
\in {\cal H}$.

\item  $\tau $-additivity: $(m(B)\eta ,\eta )=\sum (m(b_i)\eta ,\eta )$, $%
\forall \eta \in {\cal H}$, where $B=\sum B_i$ is any disjoint sum of
cardinality $|I|\leq \tau $.

\item  normalization: $m(E)=I$, where $I$ is the unit operator on ${\cal H}$.
\end{enumerate}

Then there is a quantum probability space $({\cal E},{\cal L},F)$ and a $%
\tau $-homomorphism $h:{\cal B}\to {\cal L}$ such that $h(E)F=F$, inducing
the distribution $m$: 
\[
m(B)=F^{*}h(B)F, 
\]
canonically determined on a minimal Hilbert space ${\cal E}$ generated by
the vectors $h(B)F\eta ,B\in {\cal B},\eta \in {\cal H}$. If $h_1$ and $h_2$
are equivalent observables on minimal spaces ${\cal E}_1$ and ${\cal E}_2$,
then there is a unitary transformation 
\[
U:{\cal E}_1\to {\cal E}_2,UF_1=F_2, 
\]
for which $h_2(B)=Uh_1(B)U^{*}$. \label{th:2.6}
\end{theorem}

\noindent{\sc Proof} First we prove the positive definiteness $%
\sum_{j,l}(m(B^jB^l)\eta _l,\eta _j)\geq 0$ of any operator-valued matrix $%
[k^{j,l}]_{j,l=1}^n,\ m(B^jB^l)=k^{j,l}$. Indeed, for any finite set $%
\{B^j\} $ of elements $B^j,\ j=1,2,...,n$ of the semiring ${\cal B}$, there
exists a finite partition $B^j=\sum_{B\in A^j}B$ of mutually disjoint
elements $B\in {\cal B}$ from an index set (the same for all $B^j$) $%
A\subset {\cal B},\ A^j\subseteq A,\forall j=1,2,...,n$. Defining $\eta
_{B,j}=\eta _j$, if $B\in A^j$, and $\eta _{B,j}=0$, if $B\notin A^j$, and
taking into account that $B^jB^l=\sum_{B\in A^j\cap A^l}B$, we obtain 
\begin{eqnarray*}
\sum_{j,l}(m(B^jB^l)\eta _j,\eta _l) &=&\sum_{j,l=1}^n\sum_{B\in A^j\cap
A^l}(m(B)\eta _{B,j},\eta _{B,l}) \\
&=&\sum_{j,l=1}^n\sum_{B\in A}(m(B)\eta _{B,j},\eta _{B,l}) \\
&=&\sum_{B\in A}(m(B)\sum_{j=1}^n\eta _{B,j},\sum_{j=1}^n\eta _{B,l})\geq 0
\end{eqnarray*}
owing to positivity~(\ref{eq:1.1}). Let ${\cal E}$ denote the Hilbert space
obtained by factorization and completion of the formal linear envelope of
the cartesian product ${\cal B}\times {\cal H}$ with scalar multiplication $%
\lambda (B,\eta )\mapsto (B,\lambda \eta )$ with respect to hermitian form 
\[
\langle \sum_j(B^j,\eta _j)|\sum_l(B^l,\eta _l)\rangle
=\sum_{j,l}(m(B^jB^l)\eta _l,\eta _j). 
\]
Denote by $|\sum_j(B^j,\eta _j)\rangle $ the corresponding equivalence
classes, and we express the above scalar product in the form $(\xi ^{\prime
},\xi )=\xi ^{*}\xi ^{\prime }$, where $\xi ^{*}=\langle \sum_j(B^j,\eta
_j)| $ is the linear functional corresponding to the vector $\xi
=|\sum_j(B^j,\eta _j)\rangle $. Linear operator $F:\eta \mapsto |(E,\eta
)\rangle $ defines an isometric embedding ${\cal H}\to {\cal E}$ due to the
equality 
\[
\langle (E,\eta )|(E,\eta )\rangle =(\eta ,\eta ). 
\]
For any $B\in {\cal B}$, we define a linear operator $h(B)$ on ${\cal E}$ by 
\[
h(B)|(B^{\prime },\eta )\rangle =|(BB^{\prime },\eta )\rangle . 
\]
The operator $h(B)$ is selfadjoint owing to commutativity of the semiring $%
{\cal B}$, and idempotent owing to $B^2=B$, $B\in {\cal B}$. Moreover, 
\[
h(E)F\eta =h(E)|(E,\eta )\rangle =|(E,\eta )\rangle =F\eta , 
\]
hence $h(E)F=F$. For any $B\in {\cal B},\eta \in {\cal H}$, 
\[
h(B)F\eta =|(B,\eta )\rangle , 
\]
hence $h(B)F\eta $ generate ${\cal E}$. Finally, 
\[
(F^{*}h(B)F\eta ,\eta )=(h(B)F\eta ,F\eta )=\langle (B,\eta )|(E,\eta
)\rangle =(m(B)\eta ,\eta ),\forall \eta , 
\]
so that $F^{*}h(B)F=m(B)$, $B\in {\cal B}$. To prove uniqueness, observe
that the mapping 
\[
U:h_1(B)F_1\eta \mapsto h_2(B)F_2\eta 
\]
is isometric on the generating vectors, and extends to a unitary
transformation $U:{\cal E}_1\to {\cal E}_2$ with the desired properties. The
proof is complete.

\bigskip

\ {\bf Acknowledgement} I am very grateful to Professor Sylvia Pulmanova for
the hospitality in Liptovski Jan and the encouragement to publish the
Chapter~2 of my untraslated book as the contributed paper. During the
conference and the discussion we found that some structures, such as
qusilogics and quasiobservables introduced in this chapter, have been
rediscovered then in the fuzzy logics and the measure theory on diposets. I
am indebted to Sylvia Pulmanova for pointing me out an omission of the axiom~%
$b-\left( b-a\right) =a$ in the definition of the quasilogics in the Russion
edition of the book and for the completion of the proof of the last Theorem.
I also wish to thank to her daughter, Anna for the marvellous translation of
this chapter from Russian into English, and the \TeX\ setting of the paper.

\end{document}